\documentclass[12pt]{article}
\usepackage{amssymb,amsmath}
\usepackage{mathrsfs}
\usepackage{graphicx}

\newtheorem{Thm}{Theorem}
\newtheorem{Lm}{Lemma}
\newtheorem{Rm}{Remark}

\begin{document}

\begin{center}
\large{\bf Maximum of Catalytic Branching Random Walk\\
with Regularly Varying Tails}
\end{center}
\vskip0,5cm
\begin{center}
Ekaterina Vl. Bulinskaya\footnote{ \emph{Email address:} {\tt
bulinskaya@yandex.ru}}$^,$\footnote{The work is partially supported by Russian Science Foundation under grant 17-11-01173.}
\vskip0,2cm \emph{Novosibirsk State University}
\end{center}
\vskip1cm

\begin{abstract}

For a continuous-time catalytic branching random walk (CBRW) on $\mathbb{Z}$, with an arbitrary finite number
of catalysts, we study the asymptotic behavior of position of the rightmost particle when time tends to infinity.
The mild requirements include the regular variation of the jump distribution tail for underlying random walk and
the well-known $L\log L$ condition for the offspring numbers. In our classification, given in the previous paper,
the analysis refers to supercritical CBRW. The principle result demonstrates that, after a proper normalization,
the maximum of CBRW converges in distribution to a non-trivial law. An explicit formula is provided for this
normalization and non-linear integral equations are obtained to determine the limiting  distribution function.
The novelty consists in establishing the weak convergence for CBRW with ``heavy'' tails, in contrast to the known
behavior in case of ``light'' tails of the random walk jumps. The new tools such as ``many-to-few lemma'' and spinal
decomposition appear non-efficient here. The approach developed in the paper combines the techniques of renewal theory,
Laplace transform, non-linear integral equations and large deviations theory for random sums of random variables.

\vskip0,5cm {\it Keywords and phrases}: catalytic branching random walk,
heavy tails, regular varying tails, spread of population,
$L\log L$ condition.

\vskip0,5cm 2010 {\it AMS classification}: 60J80, 60F05.

\end{abstract}

\section{Introduction}
More than fifty years ago the new stochastic models were introduced to describe simultaneously the branching of particles and their movement in space (see, e.g., a recent survey \cite{Shi_LNM_15}). Nevertheless, the variety of the model settings, imposed conditions, characteristics under consideration and relations to natural sciences explain the interest and activity in investigations of branching random walk (BRW) till nowadays. The paramount problems for BRW models are the survival of population and the rate of its spatial propagation under non-extinction condition. As for classical branching processes (see, e.g., \cite{Sev_74}), the analysis of (particles) population survival leads to the corresponding classification of models, otherwise to introduction of different regimes. Recall that, for supercritical BRW, one  observes an exponential growth of population with positive probability, contrary to critical and subcritical regimes where the population does not increase with probability one. A large part of analysis concentrates on the population growth, i.e. one considers the supercritical BRW. The rate of the population propagation in this case depends essentially on the conditions imposed on ``walking'' and ``branching''. Recent study of maximum of particles positions in the standard \emph{space-homogeneous} BRW on real line was carried out in \cite{Gantert_Hofelsauer_18}, \cite{Lifshits_14} and \cite{Mallein_16} under condition of ``light'' tails of the distribution of the random walk jump. The case of ``heavy'' tails required other handling provided, e.g., in \cite{Bhattacharya_etal_18}, \cite{Durrett_83} and \cite{GMV_2017}.

\emph{Space-inhomogeneous} models, called \emph{catalytic} (CBRW), arise when the BRW is accompanied by a set of catalysts. Now
a particle can give, with a specified probability, the offspring only at the certain points (catalysts) located in the space where a random walk takes place. With complementary probability a particle leaves the catalyst according to the mechanism of the underlying random walk. Outside the catalysts the particles perform an ordinary random walk without branching. If a particle gives (a random number of) offspring it immediately dies and the offsprings evolve in the same manner as a parent, independently of each other and of all existing particles. More details are given in Section~\ref{s:main_results}. We refer also to the papers \cite{B_JTP_14}, \cite{VTY}, \cite{Y_MCAP_17} devoted to CBRW with a finite set of catalysts and \cite{Platonova_Ryadovkin_17} dealing with a periodical infinite catalysts set. Notably the main tool in \cite{Platonova_Ryadovkin_17} and \cite{Y_MCAP_17} is the operator theory applied to the evolutionary operator related to the mean local and total particles numbers whereas the key-technique in \cite{B_JTP_14} and \cite{VTY} are probabilistic methods combining renewal theory, auxiliary branching processes, hitting times under taboo and others.

The main subject below is the rate of population propagation for continuous-time CBRW on  $\mathbb{Z}$, as time tends to infinity. This problem is solved in \cite{Carmona_Hu_14} in the case of ``light'' tails  of jump distribution of the random walk, i.e. under the Cram\'{e}r condition for jump law. Recently these investigations were extended to CBRW on $\mathbb{Z}^d$, $d\in\mathbb{N}$, in \cite{B_SPA_18}, also under the Cram\'{e}r condition. There were applied the new classification of CBRW introduced in \cite{B_TPA_15} and limit theorems for local and total particles numbers established in \cite{B_Doklady_15}. It is worth  mentioning that  related results in terms of boundedness of some random variables moments, rather than strong or weak convergence of maximum of CBRW, were obtained in a series of papers, see, e.g., \cite{Molchanov_Yarovaya_12}. Furthermore, various characteristics of particles propagation were analyzed in the framework of spatially continuous counterpart of CBRW called \emph{catalytic branching Brownian motion}, see, e.g., \cite{Bocharov_Harris_16} and \cite{Wang_Zong_17}. Although there are recent advances in the study of CBRW where the Cram\'{e}r condition fails (see, e.g., \cite{Y_MCAP_17} and references therein), there have been yet no results on convergence of maximum of CBRW under ``heavy'' tails assumption on jump distribution of the random walk. The aim of our work is to propose an approach comprising the case of mentioned ``heavy'' tails.

In fact, we find the proper normalizing factor for the maximum of CBRW on $\mathbb{Z}$ to have a non-trivial limit in a weak sense,
as time grows to infinity. Also, in the particular case of the starting point of CBRW belonging to the catalysts set, we demonstrate
that the limit cumulative distribution functions (c.d.f.) of the transformed maximum obey a finite system of non-linear integral equations. The general case of arbitrary starting point is reduced to that particular case. Having derived this system, we show that its solution exists and is unique in the relevant class of vector-functions. Within this study the right tail of jump distribution of the random walk is presupposed regularly varying and, for the distribution of the offspring number of each particle, the $L\log L$ condition holds true.

Note that our results show that the maximum of CBRW grows exponentially, in accord with the features of space-homogeneous
(without catalysts) BRW in \cite{Durrett_83} and \cite{GMV_2017}. Moreover, the restrictions on the distribution tails of
the random walk are milder than in \cite{GMV_2017}. Especially note that we derive a system of non-linear integral equations for the c.d.f. of the transformed maximum of CBRW under different starting points from the catalysts set. Afterwards we analyze the solution limit, invoking the renewal theory, Laplace transform and ``heavy'' tails theory including large deviations for random sums of random variables (see, e.g., monographs \cite{Borovkov_Borovkov_08}, \cite{FKZ_11} and \cite{Resnick_06}), etc. Surprisingly, the new technique of ``many-to-few'' lemma and spinal decomposition, helpful for the study of CBRW with ``light'' tails in \cite{B_SPA_18} and \cite{Carmona_Hu_14}, appears to be inefficient for the case of ``heavy'' tails. However, the ideas of papers \cite{Athreya_1969} and \cite{Kaplan_75}, non-related to the discussed problems, turn out to be fruitful in our investigation.

The  paper is organized as follows. In Section~\ref{s:main_results} we introduce necessary notation, formulate the main result and compare it with the previously known ones. More precisely, Theorem~\ref{T:main_result} establishes convergence in law of the normalized maximum of CBRW on $\mathbb{Z}$ and contains information on the limiting distribution. In Section~\ref{s:proof} we prove this theorem partitioning the arguments into 8 lemmas. Firstly we consider the case of a single catalyst and then extend the obtained results to the case of an arbitrary finite number of catalysts.

\section{Main Result}
\label{s:main_results}

Let us recall the description of CBRW on $\mathbb{Z}$. At the initial time ${t=0}$ there is a single particle that moves on $\mathbb{Z}$ according to a continuous-time Markov chain $S=\{S(t),t\geq0\}$ generated by the infinitesimal matrix ${Q=(q(x,y))_{x,y\in\mathbb{Z}}}$. When this particle hits a finite set of catalysts $W=\{w_1,\ldots,w_N\}\subset\mathbb{Z}$, say at the point $w_k$, it spends there random time, distributed exponentially with parameter $\beta_k>0$. Afterwards the particle either branches or leaves the point $w_k$ with probabilities $\alpha_k$ and $1-\alpha_k$ ($0\leq\alpha_k<1$), respectively. If the particle branches (at the point $w_k$), it produces a random non-negative integer number $\xi_{k}$ of offsprings, located at the same point $w_k$, and dies instantly. Whenever the particle leaves $w_k$, it jumps to the point $y\neq w_k$ with probability $-(1-\alpha_k)q(w_k,y)q(w_k,w_k)^{-1}$ and resumes its motion governed by the Markov chain $S$. All the newly born particles are supposed to behave as independent copies of their parent.

We assume that the Markov chain $S$ is irreducible and space-homogeneous, with the matrix $Q$ being conservative, i.e.,
\begin{equation}\label{condition1}
q(x,y)=q(x-y,0)=q(0,y-x)\quad\mbox{and}\quad\sum\limits_{y\in\mathbb{Z}}{q(x,y)}=0,
\end{equation}
where $q(x,y)\geq0$ for $x\neq y$ and $q:=-q(x,x)\in(0,\infty)$, for any $x,y\in\mathbb{Z}$. Stress that, contrary to, e.g., \cite{Platonova_Ryadovkin_17} and \cite{Y_MCAP_17}, we do not restrict ourselves to the case of symmetric generator $Q$. Denote by $f_k(s):={\sf E}{s^{\xi_k}}$, $s\in[0,1]$, the probability generating function of $\xi_k$, $k=1,\ldots,N$. We employ the standard assumption of existence of a finite derivative $f_k'(1)$, that is the finiteness of $m_k:={\sf E}{\xi_k}$, for any $k=1,\ldots,N$. Moreover, the $L\log L$ condition is stipulated for the offspring numbers, i.e.
\begin{equation}\label{condition2}
{\sf E}\,\xi_k\,\ln{\xi_k}<\infty,\quad k=1,\ldots,N,
\end{equation}
where, as usual, $s\,\ln{s}$ for $s=0$ equals $0$.

To formulate the main result of the paper let us introduce additional notation. The index $x$ in expressions of the form ${\sf E}_x$ and ${\sf P}_x$ marks the starting point of either CBRW or the random walk $S$, depending on the context. We temporarily forget about the catalysts and consider only the motion of a particle on $\mathbb{Z}$ according to the Markov chain $S$ with the generator $Q$ and the starting point $x$. The conditions imposed on the elements $q(x,y)$, $x,y\in\mathbb{Z}$, allow us to use an explicit construction of the random walk on $\mathbb{Z}$ with generator $Q$ (see, e.g., Theorem~1.2 in \cite{Bremaud_99}, Ch.~9, Sec.~1). Whence $S$ is a regular jump process with right continuous trajectories and, for transition times of the process, $\tau^{(0)}:=0$ and $\tau^{(n)}:=\inf\left\{t\geq\tau^{(n-1)}:S(t)\neq S(\tau^{(n-1)})\right\}$, $n\ge1$, the following property is valid. The random variables $\left\{\tau^{(n+1)}-\tau^{(n)}\right\}_{n=0}^{\infty}$ are independent and each of them has exponential distribution with parameter $q$. Denote by $\Pi=\{\Pi(t),t\geq0\}$ the Poisson process constructed as the renewal process with the interarrival times $\tau^{(n+1)}-\tau^{(n)}$, $n\in\mathbb{Z}_+$, (see, e.g., \cite{Feller_71}, Ch.~1, Sec.~4), that is, $\Pi$ is a
Poisson process with constant intensity $q$. Let $Y^i$ be the value of the $i$th jump of the random walk $S$ ($i=1,2,\ldots$). In
view of Theorem~1.2 in \cite{Bremaud_99}, Ch.~9, Sec.~1, the random variables $Y^1,Y^2,\ldots$ are i.i.d., have distribution ${\sf P}(Y^1=y)=q(0,y)/q$, $y\in\mathbb{Z}$, $y\neq0$, and do not depend on the sequence $\{\tau^{(n+1)}-\tau^{(n)}\}_{n=0}^{\infty}$. In other words, the formula
$$S(t)=x+\sum_{i=1}^{\Pi(t)}Y^i$$
is true (as usual, $\sum_{i\in\varnothing}Y^i=0$), where $x$ is the initial state of the Markov chain $S$. Due to this equality it is not difficult to show that $S$ is a process with independent increments. In what follows we consider the version of the process $S$ constructed in such a way and also often called compound Poisson process.

Set
$$\tau_x:=\mathbb{I}(S(0)=x)\inf\{t\geq0:S(t)\neq x\},$$
i.e. the stopping time $\tau_x$ is the time of the first exit from the starting point $x$ of the random walk. As usual,
$\mathbb{I}(A)$ stands for the indicator of a set $A$. Clearly, ${\sf P}_x(\tau_x\leq t)=1-e^{-qt}=:G_0(t)$, $x\in\mathbb{Z}$, $t\geq0$. Let
$$_T\overline{\tau}_{x,y}:=\mathbb{I}(S(0)=x)\inf\{t\geq0:S(t+\tau_x)=y,S(u)\notin T,\tau_x\leq u< t+\tau_x\}$$
be the time elapsed from the exit moment of this Markov chain (in other terms, particle) from the starting point
$x$ till the moment of the first hitting point $y$, whenever the particle trajectory does not pass the set $T\subset\mathbb{Z}$. Otherwise, we put ${_T\overline{\tau}_{x,y}=\infty}$. An extended random variable $_T\overline{\tau}_{x,y}$ is called
\emph{hitting time} of state $y$ \emph{under taboo} on set $T$ after exit out of starting state $x$ (see, e.g.,
\cite{B_SPL_14}). Denote by $_T\overline{F}_{x,y}(t)$, $t\geq0$, the improper c.d.f. of this extended random variable and let
$_T\overline{F}_{x,y}(\infty):=\lim_{t\to\infty}{_T\overline{F}_{x,y}(t)}$. If the taboo set $T$ is empty, expressions
$_{\varnothing}\overline{\tau}_{x,y}$ and $_{\varnothing}\overline{F}_{x,y}$ are shortened as $\overline{\tau}_{x,y}$ and $\overline{F}_{x,y}$. Mainly we will be interested in the situation when $T=W_k$, where $W_k:=W\setminus\{w_k\}$, $k=1,\ldots,N$.

Hereinafter,
$$F^{\ast}(\lambda):=\int\nolimits_{0-}^{\infty}{e^{-\lambda t}\,d{F(t)}},\quad\lambda\geq0,$$
stands for the Laplace transform of a c.d.f. $F(t)$, $t\geq0$, with support on non-negative semi-axis. For $j,k=1,\ldots,N$, $x,y\in\mathbb{Z}$ and $t\geq0$, set
$$G_j(t):=1-e^{-\beta_j t},\quad G_{j,k}(t):=G_j\ast{_{W_k}\overline{F}_{w_j,w_k}(t)},\quad{_T F_{x,y}(t)}:=G_0\ast{_T\overline{F}_{x,y}(t)},$$
where $\ast$ denotes a convolution of c.d.f. Note that by definition the function ${_T F_{x,y}(\cdot)}$ is a c.d.f. of the variable $_T\tau_{x,y}:=\tau_x+{_T\overline{\tau}_{x,y}}$ called \emph{hitting time} of state $y$ \emph{under taboo} on set $T$ when the starting state is $x$.

Consider a matrix function $D(\lambda)=\left(d_{i,j}(\lambda)\right)_{i,j=1}^N$, $\lambda\geq0$, taking values in the set of irreducible matrices of size $N\times N$, with elements defined by way of (see \cite{B_TPA_15})
$$d_{i,j}(\lambda)=\delta_{i,j}\alpha_im_iG^{\ast}_i(\lambda)+(1-\alpha_i)G^{\ast}_i(\lambda)
{_{W_j}\overline{F}^{\ast}_{w_i,w_j}(\lambda)},$$ where $\delta_{i,j}$ is the Kronecker delta. According to Definition~$1$ in
\cite{B_TPA_15} CBRW is called {\it supercritical} if the Perron root (i.e. positive eigenvalue being the spectral radius) $\rho(D(0))$ of the matrix $D(0)$ is greater than $1$. Then in view of monotonicity of
all elements of the matrix function $D(\cdot)$ there exists the solution $\nu>0$ of the equation $\rho(D(\lambda))=1$. As Theorem $1$ in \cite{B_TPA_15} shows, just this positive number $\nu$ specifies the rate of exponential growth of the mean total and local particles numbers (in the literature devoted to population dynamics and classical branching processes one traditionally speaks of Malthusian parameter). In the sequel we consider the supercritical CBRW on $\mathbb{Z}$.

Let $N(t)\subset\mathbb{Z}$ be the (random) set of particles existing in CBRW at time ${t\geq0}$. For a particle $v\in N(t)$, denote by $X_v(t)$ its position at time $t$. We study the asymptotic behavior, as $t\to\infty$, of the rightmost particle in CBRW on $\mathbb{Z}$, i.e. of the maximum of CBRW defined by $M_t:=\max\{X_v(t):v\in N(t)\}$. Naturally, set $M_t:=-\infty$ if $N(t)=\varnothing$. Our main result (see Theorem \ref{T:main_result}) shows that the normalized maximum $M_t$ converges in distribution to a non-trivial law, as $t\to\infty$, and the normalizing factor depends on the decay rate of the tails of the random walk $S$. Suppose that
\begin{equation}\label{assumption:tails}
{\sf P}(Y^1>y)=y^{-\gamma}L_1(y)=:R(y),
\end{equation}
where ${\sf P}(Y^1>y)=q^{-1}\sum\nolimits_{x=y+1}^{\infty}{q(0,x)}$ and $L_1(y)$, $y\in\mathbb{Z}$, is a slowly varying function, i.e. let the right tail of the random walk be a regularly varying function with index $-\gamma$ and $\gamma\in(0,\infty)$. Then in accordance with \cite{Seneta_76}, Ch.~1, Sec.~5, property $5^\circ$, there exists an asymptotically uniquely determined inverse function $R^{inv}(s)$, $s\geq0$, in the sense that $1/R(R^{inv}(y))\sim y$, $R^{inv}(1/R(y))\sim y$, as $y\to\infty$, $y\in\mathbb{Z_+}$, and $R^{inv}(s)=s^{1/\gamma}L_2(s)$, where $L_2(s)$, $s\geq0$, is a slowly varying function. In other words, $R^{inv}$ is an asymptotically inverse function to $1/R$. Thus, the mentioned normalizing factor for $M_t$ is chosen to be
$$L_t=R^{inv}\left(e^{\nu t}\right)=e^{\nu t/\gamma}L_2\left(e^{\nu t}\right),\quad t\geq0.$$
Moreover, we assume that, for any positive constant $c$, ``the principle of a single big jump'' (see, e.g., \cite{FKZ_11}, p.~40) is valid, i.e.
\begin{equation}\label{assumptions:tails}
{\sf P}_0(S(t)\geq c L_t)\sim t{\sf P}(Y_1\geq c L_t),\quad t\to\infty.
\end{equation}
Broad sufficient conditions for its validity are gathered, e.g., in Theorem~15.2.1 in \cite{Borovkov_Borovkov_08}.

Introduce the following function classes
$$\mathcal{C}:=\left\{\left(\varphi(\,\cdot\,;w_1),\ldots,\varphi(\,\cdot\,;w_N)\right):\varphi(\,\cdot\,;w_i)\;
\mbox{maps}\;[0,\infty)\;\mbox{onto}\;(0,1],\phantom{\frac{1}{2}}\right.$$
$$\left.\varphi(0;w_i)=1\;\mbox{and}\;\lim_{\lambda\to0+}\frac{1-\varphi(\lambda;w_i)}{\lambda}
>0,\;i=1,\ldots,N\right\},$$
$$\mathcal{C}_{\theta}:=\left\{\left(\varphi(\,\cdot\,;w_1),\ldots,\varphi(\,\cdot\,;w_N)\right)\in\mathcal{C}:
\lim_{\lambda\to0+}\frac{1-\varphi(\lambda;w_i)}{\lambda}=\theta_i,\;i=1,\ldots,N\right\}\!,$$
where $\theta=\left(\theta_1,\ldots,\theta_N\right)$, $\theta_i>0$, $i=1,\ldots,N$.

Recall the definition of the local extinction probability $Q(x,y)\!:={\sf P}_x\left(\limsup_{t\to\infty}\mu(t;y)=0\right)$, $x,y\in\mathbb{Z}$, where $\mu(t;y)$ is the number of particles in CBRW at point $y$ at time $t$ (local particles number). Theorem~2 in \cite{B_Doklady_15} asserts that the function $Q(x,y)$ depends on $x$ only, i.e. $Q(x,y)=Q(x)$, and satisfies some system of algebraic equations provided there.

\begin{Thm}\label{T:main_result}
Let, for supercritical CBRW on $\mathbb{Z}$ with Malthusian parameter $\nu$, the conditions (\ref{condition1}), (\ref{condition2}), (\ref{assumption:tails}) and (\ref{assumptions:tails}) be satisfied. Then there exists a function $\varphi(\lambda;x)$, $\lambda\geq0$, $x\in\mathbb{Z}$, such that, for each $\lambda>0$ and $x\in\mathbb{Z}$, one has
$${\sf P}_x\left(M_t/L_t\leq\lambda^{-1/\gamma}\right)\to\varphi(\lambda;x)\in\left(0,1\right),\quad t\to\infty,$$
where $\varphi(\lambda;x)\to Q(x)$, as $\lambda\to+\infty$, for each $x\in\mathbb{Z}$. Moreover, for $x\in\mathbb{Z}\setminus W$,
the function $\varphi(\lambda;x)$, $\lambda\geq0$, admits the following representation
$$\varphi(\lambda;x)=\sum_{k=1}^N{\int\nolimits_0^{\infty}{\varphi(\lambda e^{-\nu u};w_k)\,d {_{W_k}F_{x,w_k}(u)}}}+1-\sum_{k=1}^N{_{W_k}F_{x,w_k}(\infty)},$$
where functions $\varphi(\,\cdot\,;w_j)$, $j=1,\ldots,N$, satisfy the system of integral equations
\begin{eqnarray}\label{varphi(lambda,wj)=system_equations}
& &\varphi(\lambda;w_j)=\alpha_j\int\nolimits_0^{\infty}{f_j(\varphi(\lambda e^{-\nu u};w_j))\,dG_j(u)}\\
&+&(1-\alpha_j)\sum_{k=1}^N{\int\nolimits_0^{\infty}{\varphi(\lambda e^{-\nu u};w_k)\,dG_{j,k}(u)}}
+(1-\alpha_j)\!\left(\!1\!-\!\sum_{k=1}^N{_{W_k}F_{w_j,w_k}(\infty)}\!\right)\!.\nonumber
\end{eqnarray}
The equations system (\ref{varphi(lambda,wj)=system_equations}) has a unique solution in the function class $\mathcal{C}_{\theta}$, for each $\theta=\left(\theta_1,\ldots,\theta_N\right)$, $\theta_i>0$, $i=1,\ldots,N$.
\end{Thm}

The function $\varphi$ has already emerged in Theorem~4 in \cite{B_Doklady_15} as the Laplace transform of the limiting distribution arising in the limit theorem for total and local particles numbers in CBRW. Thus, to obtain a known limit  in Theorem~\ref{T:main_result} we consider a variable $\lambda^{-1/\gamma}$ instead of $\lambda$ (and $M_t/L_t$ instead of $(M_t/L_t)^{-\gamma}$, since the latter expression is undefined when $M_t<0$).

Moreover, as stated in Theorem~4 of \cite{B_Doklady_15}, both total and local particles numbers, being normalized by their means, converge in distribution to a non-degenerate random variable which vanishes with probability $Q(x)$ of \emph{local extinction} of population in CBRW starting at the point $x$. Similarly, as shows Theorem~\ref{T:main_result}, the trivial relation $M_t/L_t\to0$, as $t\to\infty$, is only realized with the same probability $Q(x)$ of \emph{local extinction} of population in CBRW with starting point $x$. This means that the normalizing factor $L_t$ is determined aptly.

Note that further Lemma~\ref{L:limlim=theta} and its counterpart in case of multiple catalysts provide the value $\lim_{\lambda\to0+}\lim_{t\to\infty}\lambda^{-1}{\sf P}_{w_j}\left(M_t/L_t>\lambda^{-1/\gamma}\right)$ which we denote by $\theta_j>0$, $j=1,\ldots,N$, with $\theta=\left(\theta_1,\ldots,\theta_N\right)$. Therefore, according to Theorem~\ref{T:main_result} the limit of ${\sf P}_{w_j}\left(M_t/L_t\leq\lambda^{-1/\gamma}\right)$, $j=1,\ldots,N$, as $t\to\infty$, is uniquely determined as the solution $\left(\varphi(\cdot\,;w_1),\ldots,\varphi(\cdot\,;w_N)\right)$ to the system (\ref{varphi(lambda,wj)=system_equations}) in the class $\mathcal{C}_{\theta}$.

\begin{Rm}
One should compare the results of paper \cite{Carmona_Hu_14} and Theorem~\ref{T:main_result} treating the behavior of maximum in case of ``light'' and ``heavy'' tails, respectively. Firstly, the normalizing factors are different, linear versus exponential. Secondly,  \cite{Carmona_Hu_14} establishes almost sure convergence of $M_t/t$ to a constant and we show that distribution of  $M_t/L_t$ weakly converges to a non-degenerate law. Note also that, in contrast to  \cite{Carmona_Hu_14}, we do not assume even the existence of the expectation of random walk jumps. Recall that for i.i.d. random summands the existence of expectation implies strong LLN, whereas when expectation lacks but the summands have regularly varying tails, one can only ensure weak convergence to a stable law. Apparently the same effect might be behind the {\bf strong} convergence of normalized ``maximum'' stated in \cite{B_SPA_18} and \cite{Carmona_Hu_14} versus the present {\bf weak} convergence result under other normalization.
\end{Rm}

\section{Proof of Theorem~\ref{T:main_result}}\label{s:proof}

\emph{By the subsequent Lemmas~\ref{L:equation}--\ref{L:limK=0}, we establish Theorem~\ref{T:main_result} in case of a single catalyst $w_1$ located, without loss of generality, at the origin, that is $W=\{w_1\}$ with $w_1=0$, and the starting point of CBRW being $0$ as well. Later we will turn to the general case.}

\vskip0.1cm

The first lemma provides an integral equation for the tail of c.d.f. of the maximum $M_t$ of CBRW.

\begin{Lm}\label{L:equation}
Let condition (\ref{condition1}) be valid. Then the probability ${\sf P}_0(M_t>u)$, $t\geq0$, $u\in\mathbb{R}$, satisfies the following non-linear integral equation of convolution type
\begin{eqnarray}
{\sf P}_0(M_t>u)&=&\alpha_1\int\nolimits_0^t{\left(1-f_1\left(1-{\sf P}_0(M_{t-s}> u)\right)\right)\,dG_1(s)}\label{P_0(Mt_>_Lt_u)_equation}\\
&+&(1-\alpha_1)\int\nolimits_0^t{{\sf P}_0(M_{t-s}>u)\,dG_{1,1}(s)}+I(t;u),\nonumber
\end{eqnarray}
where, for $u\in\mathbb{R}$, one has
$$I(t;u):=\sum_{y\neq0}{(1-\alpha_1)\frac{q(0,y)}{q}\int\nolimits_0^t{{\sf P}_y\left(S(t-s)>u,\tau_{y,0}>t-s\right)\,dG_1(s)}}.$$
\end{Lm}
{\sc Proof.}
Consider all the possible evolutions of the parent particle in CBRW. Namely, after time, distributed exponentially with parameter $\beta_1$, it may either produce $k\in\mathbb{Z}_+$ offsprings with probability $\alpha_1{\sf P}(\xi_1=k)$, or jump to the point $y\neq 0$ with probability $(1-\alpha_1)q(0,y)/q$ and afterwards first return to the origin at time $\overline{\tau}_{0,0}$. If the parent particle does not return to the origin until time $t$, it performs an ordinary random walk $S$ starting from $y$. At last, it might occur that the parent particle has not undergone changes by time $t$. Summarizing all the above we can write the following formula, for any $u\in\mathbb{R}$,
\begin{eqnarray*}
{\sf P}_0(M_t\leq u)&=&\alpha_1\sum_{k=0}^{\infty}{\sf P}(\xi_1=k)\int\nolimits_0^t{\left({\sf P}_0(M_{t-s}\leq u)\right)^k\,dG_1(s)}+(1-G_1(t))\\
&+&\sum_{y\neq 0}{(1-\alpha_1)\frac{q(0,y)}{q}\int\nolimits_0^t{{\sf P}_0(M_{t-s}\leq u)\,d\left(G_1\ast\overline{F}_{0,0}(s)\right)}}\\
&+&\sum_{y\neq 0}{(1-\alpha_1)\frac{q(0,y)}{q}\int\nolimits_0^t{{\sf P}_y\left(S(t-s)\leq u,\,\tau_{y,0}>t-s\right)\,dG_1(s)}}.
\end{eqnarray*}
Rewriting the latter equation with respect to unknown function ${\sf P}_0\left(M_t>u\right)$ and taking into account the obvious identity
$${\sf P}_y\left(S(s)\leq u,\tau_{y,0}>s\right)=1-F_{y,0}(s)-{\sf P}_y\left(S(s)>u,\tau_{y,0}>s\right),\quad s\geq0,$$
we get the assertion of Lemma~\ref{L:equation}. $\square$

The following lemma provides a convenient form for the function $I$ expressed in terms of the probability ${\sf P}_0\left(S(t)>u\right)$ when $u\geq0$.

\begin{Lm}\label{L:J-1(t;a)=}
Let condition (\ref{condition1}) be satisfied. Then, for any $t,u\geq0$, the following identity holds true
\begin{eqnarray}\label{J-1(t;a)=}
\frac{q I(t;u)}{(1-\alpha_1)\beta_1}&=&{\sf P}_0\left(S(t)>u\right)-\int\nolimits_0^t{{\sf P}_0\left(S(t-s)>u\right)\,d F_{0,0}(s)}\\
&-&\frac{\beta_1-q}{\beta_1}\int\nolimits_0^t{{\sf P}_0\left(S(t-s)>u\right)\,d\left(G_1(s)-G_1\ast F_{0,0}(s)\right)}.\nonumber
\end{eqnarray}
\end{Lm}
{\sc Proof.}
It is not difficult to see that, for any non-negative $u$, one has
\begin{equation}\label{identity_for_random_walk}
\sum_{y\neq 0}{\frac{q(0,y)}{q}\int\nolimits_0^t{{\sf P}_y\left(S(t-s)>u,\,\tau_{y,0}>t-s\right)\,dG_0(s)}}={\sf P}_0\left(S(t)>u,\,\tau_{0,0}>t\right).
\end{equation}
Observe that the function $I$ takes form of the left-hand side of the latter identity upon omitting the factor $(1-\alpha_1)$ and replacing the function $G_1$ by $G_0$ in the definition of $I$. On this footing let us first verify that $(1-\alpha_1)^{-1}\beta_1^{-1}qI(t;u)$, for $u\geq0$, equals the following expression
\begin{equation}\label{J-1=alternative-representation}
{\sf P}_0\left(S(t)>u,\tau_{0,0}>t\right)-\frac{\beta_1-q}{\beta_1}\int\nolimits_0^t{{\sf P}_0\left(S(t-s)>u,\tau_{0,0}>t-s\right)dG_1(s)}.
\end{equation}
Indeed, implementing the Laplace transform of the latter expression, one can write
$$\frac{\lambda+q}{\lambda+\beta_1}\int\nolimits_0^{\infty}{e^{-\lambda t}{\sf P}_0\left(S(t)>u,\tau_{0,0}>t\right)\,dt},$$
since $G_1^{\ast}(\lambda)=\beta_1/\left(\lambda+\beta_1\right)$, $\lambda\geq0$.
We obtain the same expression by applying the Laplace transform to the function $(1-\alpha_1)^{-1}\beta_1^{-1}qI(t;u)$ and also taking into account identity (\ref{identity_for_random_walk}). Further on, the Laplace transform uniqueness entails the alternative representation (\ref{J-1=alternative-representation}) for $(1-\alpha_1)^{-1}\beta_1^{-1}qI(t;u)$.

Evidently, for any $u\geq0$ one has
$${\sf P}_0\left(S(t)>u,\tau_{0,0}>t\right)={\sf P}_0\left(S(t)>u\right)-\int\nolimits_0^t{{\sf P}_0\left(S(t-s)>u\right)\,d F_{0,0}(s)}.$$
Substituting the latter equality into the verified alternative representation (\ref{J-1=alternative-representation}) for the function $(1-\alpha_1)^{-1}\beta_1^{-1}qI(t;u)$, we come to the assertion of Lemma~\ref{L:J-1(t;a)=}. $\square$

Recall that by the definition of supercritical regime of CBRW the relations $\alpha_1m_1+(1-\alpha_1)F_{0,0}(\infty)>1$ and $\alpha_1m_1G_1^{\ast}(\nu)+(1-\alpha_1)\,G_1^{\ast}(\nu)\overline{F}^{\,\ast}_{0,0}(\nu)=1$ are valid. In terms of the function $G(t):=\alpha_1m_1G_1(t)+(1-\alpha_1)\,G_1\ast\overline{F}_{0,0}(t)$, $t\geq0$, it means that $G^{\ast}(\nu)=1$.

\begin{Lm}\label{L:J-1(t;lambdaLt)ast_sim}
Whenever assumptions (\ref{condition1}), (\ref{assumption:tails}) and (\ref{assumptions:tails}) hold true, one has, for each $\lambda>0$ and $r\geq0$,
\begin{equation}\label{J-1(t;lambdaLt)ast_sim}
\int\nolimits_0^t{I(t-u;\lambda^{-1/\gamma}L_{t+r})\,d\sum_{k=0}^{\infty}G^{\ast k}(u)}\sim K\lambda e^{-\nu r},\quad t\to\infty,
\end{equation}
where the constant $K$ is of the following form
$$K=\frac{(1-\alpha_1)\beta_1(\nu+q)}{q(\nu+\beta_1)}\left(1-F^{\,\ast}_{0,0}(\nu)\right)\int\nolimits_0^{\infty}{se^{-\nu s}\,ds}\left(\int\nolimits_0^{\infty}{se^{-\nu s}\,dG(s)}\right)^{-1}.$$
\end{Lm}
{\sc Proof.}
Firstly, for each $\lambda>0$, $r\geq0$, let us study the asymptotic behavior of $I(t;\lambda^{-1/\gamma}L_{t+r})$ as ${t\to\infty}$. Assumption (\ref{assumptions:tails}) implies that
$${\sf P}_0\left(S(t-s)>\lambda^{-1/\gamma}L_{t+r}\right)\sim(t-s)R\left(\lambda^{-1/\gamma}L_{t+r}\right)\sim(t-s)\lambda\, e^{-\nu(t+r)},$$
as $t-s\to\infty$, whenever $\lambda>0$ and $r\geq0$. Therefore, identity (\ref{J-1(t;a)=}) yields
\begin{eqnarray}\label{I(t,lambdaL_t)sim}
& &(1-\alpha_1)^{-1}\beta_1^{-1}q I(t;\lambda^{-1/\gamma}L_{t+r})\sim\lambda\,e^{-\nu(t+r)}\\
&\times&\left(t-\int\nolimits_0^t{(t-s)\,d F_{0,0}(s)}-\frac{\beta_1-q}{\beta_1}\int\nolimits_0^t{(t-s)\,d\left(G_1(s)-G_1\ast F_{0,0}(s)\right)}\right),\nonumber
\end{eqnarray}
as $t\to\infty$. In view of Theorem~25 in \cite{Vatutin_book_09}, p.~30, it follows that
$$\int\nolimits_0^t{I(t-u;\lambda^{-1/\gamma}L_{t+r})\,d\sum_{k=0}^{\infty}G^{\ast k}(u)}\sim\frac{(1-\alpha_1)\beta_1\lambda\,e^{-\nu(t+r)}e^{\nu t}}{q\int\nolimits_0^{\infty}{se^{-\nu s}\,dG(s)}}$$
$$\times\int\limits_0^{\infty}{\!e^{-\nu s}\!\!\left(\!\!s-\!\int\limits_0^s{\!(s-u)\,d F_{0,0}(u)}\!-\!\frac{\beta_1-q}{\beta_1}\!\int\limits_0^s{\!(s-u)\,d\left(G_1(u)\!-\!G_1\ast F_{0,0}(u)\right)}\!\!\right)\!ds}$$
$$=K\lambda e^{-\nu r},$$
as $t\to\infty$. Here we employ the Laplace transform of convolutions property and the formula $\int\nolimits_0^{\infty}{e^{-\nu s}\,dG_1(s)}=\beta_1/(\nu+\beta_1)$. Lemma~\ref{L:J-1(t;lambdaLt)ast_sim} is proved completely. $\square$

Next we derive an upper bound for the probability ${\sf P}_0\left(M_t>\lambda^{-1/\gamma}L_{t+r}\right)$.

\begin{Lm}\label{L:lower_estimate}
If conditions (\ref{condition1}), (\ref{assumption:tails}) and (\ref{assumptions:tails}) are satisfied, then,
for any $r,t\geq0$, $\lambda>0$ and some positive constant $C$, the following inequality is valid
\begin{equation}\label{P_0(Mt_>_Lt+r_u)_upper_estimate}
{\sf P}_0\left(M_t>\lambda^{-1/\gamma}L_{t+r}\right)\leq C\lambda e^{-\nu r}.
\end{equation}
\end{Lm}
{\sc Proof.}
For any $u\in\mathbb{R}$, by mean value theorem on $f_1$, equation (\ref{P_0(Mt_>_Lt_u)_equation}) entails the inequality
\begin{equation}
{\sf P}_0\!\left(M_t>u\right)\!\leq\!\int\nolimits_0^t{{\sf P}_0\!\left(M_{t-s}>u\right)\,d G(s)}+I\!\left(t;u\right)\!.
\end{equation}
Iterating this inequality $k$ times we get
$${\sf P}_0(M_t>u)\leq\int\nolimits_0^t{{\sf P}_0\left(M_{t-s}>u\right)\,d G^{\ast (k+1)}(s)}+\int\nolimits_0^t{I(t-s;u)\,d\sum_{j=0}^k{G^{\ast j}(s)}}.$$
For any fixed $t$, one has $G^{\ast k}(t)\to 0$, as $k\to\infty$. For example, this is due to Lemma~22 in \cite{Vatutin_book_09}. Hence, the term $\int\nolimits_0^t{{\sf P}_0\left(M_{t-s}>u\right)\,dG^{\ast(k+1)}(s)}$ is negligibly small for large $k$. Therefore, the latter inequality can be rewritten as follows
\begin{equation}\label{P_0(Ltgamma/Mtgamma<lambda)_inequality}
{\sf P}_0(M_t>u)\leq\int\nolimits_0^t{I(t-s;u)\,d\sum_{j=0}^{\infty}{G^{\ast j}(s)}}.
\end{equation}
Letting $u=\lambda^{-1/\gamma}L_{t+r}$ in this relation and invoking Lemma~\ref{L:J-1(t;lambdaLt)ast_sim}, we come to the statement of Lemma~\ref{L:lower_estimate}. $\square$

Denote by $J\left(t;u\right)$, $t\geq0$, $u\in\mathbb{R}$, the difference
$$m_1\int\nolimits_0^t{{\sf P}_0\left(M_{t-s}>u\right)\,dG_1(s)}-\int\nolimits_0^t{\left(1-f_1\left(1-{\sf P}_0\left(M_{t-s}>u\right)\right)\right)\,dG_1(s)}.$$

\begin{Lm}\label{L:I-1(t;lambdaLt)ast_sim}
If conditions (\ref{condition1}), (\ref{condition2}), (\ref{assumption:tails}) and (\ref{assumptions:tails}) are satisfied, then the following relation holds true
$$\lim_{\lambda\to0+}\lim_{t\to\infty}\frac{1}{\lambda}\int\nolimits_0^t{J\left(t-s;\lambda^{-1/\gamma}L_t\right)\,d\sum_{j=0}^{\infty}{G^{\ast j}(s)}}=0.$$
\end{Lm}
{\sc Proof.} Mean value theorem on $f_1$ and Lemma~\ref{L:lower_estimate}, applied when $C\lambda<1$, ensure that
$$0\leq\frac{1}{\lambda}\int\nolimits_0^t{J\left(t-s;\lambda^{-1/\gamma}L_t\right)\,d\sum_{j=0}^{\infty}{G^{\ast j}(s)}}$$
$$\leq C\int\nolimits_0^t{\left(m_1-f_1'\left(1-C\lambda e^{-\nu s}\right)\right)e^{-\nu s}\,d\left(G_1\ast\sum_{j=0}^{\infty}{G^{\ast j}}(s)\right)}.$$
By the definition of the Malthusian parameter $\nu$, in view of Theorem~25 in \cite{Vatutin_book_09}, p.30, one has, as $t\to\infty$,
$$G_1\ast\sum_{j=0}^{\infty}{G^{\ast j}}(t)\sim e^{\nu t}\frac{\int_0^{\infty}{\left(1-e^{-\beta_1 s}\right)e^{-\nu s}\,ds}}{\int_0^{\infty}{se^{-\nu s}\,dG(s)}}.$$
Returning to the previous chain of inequalities, we see that
$$\frac{1}{\lambda}\int\nolimits_0^t{J(t-s;\lambda^{-1/\gamma}L_t)\,d\sum_{j=0}^{\infty}{G^{\ast j}(s)}}\leq C_1\int\nolimits_0^t{\left(m_1-f_1'\left(1-C\lambda e^{-\nu s}\right)\right)\,ds},$$
for some positive constant $C_1$.
Let us show that the latter integral converges, as $t\to\infty$, whenever ${\sf E}\,\xi_1\ln{\xi_1}<\infty$. Indeed,
$$\int\nolimits_0^t{\left(m_1-f_1'\left(1-C\lambda e^{-\nu s}\right)\right)\,ds}=\int\nolimits_0^t{\left({\sf E}\xi_1-{\sf E}\left(\xi_1\left(1-C\lambda e^{-\nu s}\right)^{\xi_1-1}\right)\right)\,ds}$$
$$=\nu^{-1}{\sf E}\left(\xi_1\int\nolimits_{1-C\lambda}^{1-C\lambda e^{-\nu t}}{\frac{1-u^{\xi_1-1}}{1-u}\,du}\right)=\nu^{-1}{\sf E}\left(\xi_1\int\nolimits_{1-C\lambda}^{1-C\lambda e^{-\nu t}}{\sum_{k=1}^{\xi_1-1}u^{k-1}}\,du\right)$$
$$=\nu^{-1}{\sf E}\!\left(\xi_1\!\sum_{k=1}^{\xi_1-1}\!\frac{\left(1-C\lambda e^{-\nu t}\right)^k-\left(1-C\lambda\right)^k}{k}\right)\!\leq\nu^{-1}{\sf E}\!\left(\xi_1\!\sum_{k=1}^{\xi_1-1}\!\frac{1-\left(1-C\lambda\right)^k}{k}\!\right)$$
$$\leq\nu^{-1}{\sf E}\left(\xi_1\left(1-\left(1-C\lambda\right)^{\xi_1}\right)\sum_{k=1}^{\xi_1-1}\frac{1}{k}\right)\leq\nu^{-1}{\sf E}\left(\xi_1\ln{\xi_1}\left(1-\left(1-C\lambda\right)^{\xi_1}\right)\right).$$
Here we performed the variable change $u=1-C\lambda e^{-\nu s}$, whence $du=C\lambda\nu e^{-\nu s}\,ds$, i.e. ${ds=du/\left(\nu(1-u)\right)}$. We obtain ${\sf E}\left(\xi_1\ln{\xi_1}\left(1-\left(1-C\lambda\right)^{\xi_1}\right)\right)\to0$, as $\lambda\to0+$, applying the bounded convergence theorem. Thus, Lemma~\ref{L:I-1(t;lambdaLt)ast_sim} is proved completely. $\square$

Lemma~\ref{L:lower_estimate} implies $\lim_{\lambda\to0+}{\sf P}_0\left(M_t>\lambda^{-1/\gamma}L_{t+r}\right)=0$ and $\lambda^{-1}{{\sf P}_0\left(M_t>\lambda^{-1/\gamma}L_{t+r}\right)\leq C e^{-\nu r}}$, $\lambda>0$, $r,t\geq0$. The next lemma refines the latter assertion if $t\to\infty$ and then $\lambda\to0+$.

\begin{Lm}\label{L:limlim=theta}
If conditions (\ref{condition1}), (\ref{condition2}), (\ref{assumption:tails}) and (\ref{assumptions:tails}) hold true, then the following relation is valid
\begin{equation}
\lim_{\lambda\to0+}\lim_{t\to\infty}\frac{{\sf P}_0\left(M_t>\lambda^{-1/\gamma}L_t\right)}{\lambda}=K.
\end{equation}
\end{Lm}
{\sc Proof.}
In view of (\ref{P_0(Mt_>_Lt_u)_equation}), for any $u\in\mathbb{R}$, one has
$${\sf P}_0\left(M_t>u\right)=\int\nolimits_0^t{{\sf P}_0\left(M_{t-s}>u\right)\,d G(s)}
+I(t;u)-J(t;u).$$
Iterating this equation $k$ times yields
\begin{eqnarray*}
{\sf P}_0\left(M_t>u\right)&=&\int\nolimits_0^t{{\sf P}_0\left(M_{t-s}>u\right)\,dG^{\ast(k+1)}(s)}\\
&+&\int\nolimits_0^t{I(t-s;u)\,d\sum_{j=0}^k{G^{\ast j}(s)}}-\int\nolimits_0^t{J(t-s;u)\,d\sum_{j=0}^k{G^{\ast j}(s)}}.
\end{eqnarray*}
For any fixed $t$, again by Lemma~22 in \cite{Vatutin_book_09}, one has $G^{\ast k}(t)\to 0$, as $k\to\infty$. Hence, the term $\int\nolimits_0^t{{\sf P}_0\left(M_{t-s}>u\right)\,dG^{\ast(k+1)}(s)}$ is negligibly small for large $k$. Therefore, the latter equation can be rewritten as follows
$${\sf P}_0\left(M_t>\lambda^{-1/\gamma}L_t\right)=\int\nolimits_0^t{I\left(t-s;\lambda^{-1/\gamma}L_t\right)\,d\sum_{j=0}^{\infty}{G^{\ast j}(s)}}-\int\nolimits_0^t{J\left(t-s;\lambda^{-1/\gamma}L_t\right)\,d\sum_{j=0}^{\infty}{G^{\ast j}(s)}}.$$
Dividing by $\lambda$ both parts of the derived equality, then tending $t$ to infinity and subsequently $\lambda$ to $0$ from the right, we deduce the assertion of Lemma~\ref{L:limlim=theta} in view of Lemmas~\ref{L:J-1(t;lambdaLt)ast_sim} and \ref{L:I-1(t;lambdaLt)ast_sim}. $\square$

We temporarily return to the case of arbitrary $N$ and prove the functional theory part of our main result Theorem~\ref{T:main_result}.

\begin{Lm}\label{L:integral_equation}
If conditions (\ref{condition1}) and (\ref{condition2}) are satisfied, then the equations system (\ref{varphi(lambda,wj)=system_equations}) has a unique solution in the function class $\mathcal{C}_{\theta}$, for each $\theta=\left(\theta_1,\ldots,\theta_N\right)$, $\theta_i>0$, $i=1,\ldots,N$.
\end{Lm}
{\sc Proof.}
In case of $N=1$ and the starting point $x=w_1$, the proof of Lemma~\ref{L:integral_equation} mainly repeats those of Theorems~1, 2 and 3 in \cite{Athreya_1969}, whereas, in case of multiple catalysts and $x\in W$, the proof repeats their generalizations in \cite{Kaplan_75}, Theorems~2.1, 2.2 and 2.3. While proving we essentially base on the fact that the Perron root of $D(0)$ is greater than $1$ in view of the supercritical regime under consideration. Furthermore, we heavily employ the definition of the Malthusian parameter $\nu$ and the Frobenius theory. Since main ideas of the argument justifying Lemma~\ref{L:integral_equation} are exploited below while establishing Lemma~\ref{L:limK=0}, the remaining details of Lemma~\ref{L:integral_equation} proof are omitted. $\square$

The next lemma coincides with the statement of Theorem~\ref{T:main_result} when $N=1$, $w_1=0$ and the starting point of CBRW is $x=0$.

\begin{Lm}\label{L:limK=0}
Let conditions (\ref{condition1}), (\ref{condition2}), (\ref{assumption:tails}) and (\ref{assumptions:tails}) be valid. Then, for each $\lambda>0$,
$$\lim_{t\to\infty}{\frac{{\sf P}_0\left(M_t\leq\lambda^{-1/\gamma}L_t\right)-\varphi(\lambda;0)}{\lambda}}=0,$$
where $\varphi(\cdot;0)\in\mathcal{C}_K$.
\end{Lm}
{\sc Proof.}
Let $K(t;\lambda)$ stand for $\lambda^{-1}\left(1-{\sf P}_0\left(M_t>\lambda^{-1/\gamma}L_t\right)-\varphi(\lambda;0)\right)$. Firstly note that
\begin{equation}\label{limlimK(t,lambda)=0}
\lim_{\lambda\to0+}\lim_{t\to\infty}|K(t;\lambda)|=0.
\end{equation}
This is true by virtue of Lemmas~\ref{L:limlim=theta}, \ref{L:integral_equation} and the triangle inequality
$$|K(t;\lambda)|\leq\left|\frac{{\sf P}_0\left(M_t>\lambda^{-1/\gamma}L_t\right)}{\lambda}-K\right|+\left|\frac{1-\varphi(\lambda;0)}{\lambda}-K\right|.$$

To prove the desired statement it is sufficient to verify that $K(\lambda):=\lim_{T\to\infty}K_T(\lambda)=0$, where $K_T(\lambda):=\sup_{t\geq T}|K(t;\lambda)|$. Equations (\ref{varphi(lambda,wj)=system_equations}) (when $N=1$ and $w_1=x=0$) and (\ref{P_0(Mt_>_Lt_u)_equation}) imply that
\begin{equation}\label{K(t;lambda)leqI2+I3}
\frac{{\sf P}_0\left(M_t>\lambda^{-1/\gamma}L_t\right)-1+\varphi(\lambda;0)}{\lambda}=I_{11}(t,T;\lambda)+I_{12}(t,T;\lambda)+I_2(t;\lambda),
\end{equation}
where, for $T<t$, we set $I_{11}(t,T;\lambda)$  equal to
\begin{eqnarray*}
& &\frac{\alpha_1}{\lambda}\int\nolimits_0^{t-T}{\left(f_1(\varphi(\lambda e^{-\nu s};0))-f_1\left(1-{\sf P}_0\left(M_{t-s}>\lambda^{-1/\gamma}L_t\right)\right)\right)\,dG_1(s)}\\
&+&\frac{1-\alpha_1}{\lambda}\int\nolimits_0^{t-T}{\left({\sf P}_0\left(M_{t-s}>\lambda^{-1/\gamma}L_t\right)-1+\varphi(\lambda e^{-\nu s};0)\right)\,d\left(G_1\ast\overline{F}_{0,0}(s)\right)},
\end{eqnarray*}
also $I_{12}(t,T;\lambda)$ differs from $I_{11}(t,T;\lambda)$ by the interval of integration only, i.e. $\int_{t-T}^t$ appears instead of $\int_0^{t-T}$, and finally
\begin{eqnarray*}
I_2(t;\lambda)&:=&\lambda^{-1}I\left(t;\lambda^{-1/\gamma}L_t\right)-\frac{\alpha_1}{\lambda}\int\nolimits_t^{\infty}
{\left(1-f_1\left(\varphi\left(\lambda e^{-\nu s};0\right)\right)\right)\,dG_1(s)}\\
&-&\frac{1-\alpha_1}{\lambda}\int\nolimits_t^{\infty}{\left(1-\varphi\!\left(\lambda e^{-\nu s};0\right)\right)\,d\left(G_1\ast\overline{F}_{0,0}(s)\right)}.
\end{eqnarray*}

It follows from relation (\ref{I(t,lambdaL_t)sim}) that $\lambda^{-1}I\left(t;\lambda^{-1/\gamma}L_t\right)\leq C_2 te^{-\nu t}$ for some positive constant $C_2$. Therefore, on account of mean value theorem on $f_1$ and the boundedness of function $\left(1-\varphi(\lambda;0)\right)/\lambda$, $\lambda\geq0$, by some constant $C_3\geq K$, we have
\begin{eqnarray}
& &\left|I_2(t;\lambda)\right|\leq C_2 te^{-\nu t}+\alpha_1 m_1\int\nolimits_t^{\infty}{\frac{1-\varphi(\lambda e^{-\nu s};0)}{\lambda e^{-\nu s}}e^{-\nu s}}\,d G_1(s)\label{I3leq}\\
&+&\!(1-\alpha_1)\!\int\limits_t^{\infty}{\!\frac{1-\varphi(\lambda e^{-\nu s};0)}{\lambda e^{-\nu s}}e^{-\nu s}\,d\left(G_1\ast\overline{F}_{0,0}(s)\right)}\!\leq\! C_2 te^{-\nu t}+C_3\left(1-\widetilde{G}(t)\right)\nonumber.
\end{eqnarray}
Here $\widetilde{G}$ is a c.d.f. such that $d\widetilde{G}(s)=e^{-\nu s}\,dG(s)$, $s\geq0$.

Let $t>2T$. Then, by virtue of (\ref{limlimK(t,lambda)=0}), mean value theorem on $f_1$ and Lemma~\ref{L:lower_estimate}, we obtain (for some positive constant $C_4$) the relation
\begin{eqnarray}\label{I22leq}
& &\left|I_{12}(t,T;\lambda)\right|\leq\lambda^{-1}\int\nolimits_{t-T}^t{\left|{\sf P}_0\left(M_{t-s}>\lambda^{-1/\gamma}L_t\right)-1+\varphi(\lambda e^{-\nu s};0)\right|\,d G(s)}\nonumber\\
&\leq&\lambda^{-1}\int\nolimits_{t-T}^t{\left|{\sf P}_0\left(M_{t-s}>\lambda^{-1/\gamma}L_t\right)-{\sf P}_0\left(M_{t-s}>\lambda^{-1/\gamma}e^{\nu s/\gamma}L_{t-s}\right)\right|\,d G(s)}\nonumber\\
&+&\int\nolimits_{t-T}^t{\left|K(t-s;\lambda e^{-\nu s})\right|e^{-\nu s}\,d G(s)}\leq C_4\left(1-\widetilde{G}(T)\right).
\end{eqnarray}

It follows from Lemma~\ref{L:lower_estimate} and its proof that, for any $\varepsilon>0$, there exists $t_0$ such that, for $t\geq t_0$, $s\leq t-T$, one has
$$\left|{\sf P}_0\left(M_{t-s}>\lambda^{-1/\gamma}L_t\right)-{\sf P}_0\left(M_{t-s}>\lambda^{-1/\gamma}e^{\nu s/\gamma}L_{t-s}\right)\right|\leq\varepsilon\lambda e^{-\nu s}.$$
Hence, for any $t>T$, again, in view of mean value theorem on $f_1$, we infer that
\begin{eqnarray}\label{I12leq}
& &\left|I_{11}(t,T;\lambda)\right|\leq{\lambda}^{-1}\int\nolimits_0^{t-T}{\left|{\sf P}_0\left(M_{t-s}>\lambda^{-1/\gamma}L_t\right)-1+\varphi(\lambda e^{-\nu s};0)\right|\,d G(s)}\nonumber\\
&\leq&{\lambda}^{-1}\int\nolimits_0^{t-T}{\left|{\sf P}_0\left(M_{t-s}>\lambda^{-1/\gamma}L_t\right)-{\sf P}_0\left(M_{t-s}>\lambda^{-1/\gamma}e^{\nu s/\gamma}L_{t-s}\right)\right|\,d G(s)}\nonumber\\
&+&\int\nolimits_0^{t-T}{\left|K(t-s;\lambda e^{-\nu s})\right|e^{-\nu s}\,d G(s)}\leq\varepsilon\widetilde{G}\left(t-T\right)\nonumber\\
&+&\int_0^{t-T}{K_T(\lambda e^{-\nu s})\,d\widetilde{G}(s)}\leq\varepsilon+{\sf E}{K_T\left(\lambda e^{-\nu\zeta}\right)},
\end{eqnarray}
where $\zeta$ is a random variable with $\widetilde{G}$ as its c.d.f.

Combination of relations (\ref{K(t;lambda)leqI2+I3})--(\ref{I12leq}), for $t>2T$, leads to the inequality
$$\left|K(t;\lambda)\right|\leq C_2te^{-\nu t}+C_3\left(1-\widetilde{G}(t)\right)+C_4\left(1-\widetilde{G}(T)\right)+\varepsilon+{\sf E}{K_T\left(\lambda e^{-\nu\zeta}\right)}.$$
It means that
$$K_{2T}(\lambda)\leq{\sf E}{K_T\left(\lambda e^{-\nu\zeta}\right)}+\varepsilon+C_2Te^{-\nu T}+\left(C_3+C_4\right)\left(1-\widetilde{G}(T)\right).$$
Letting $T\to\infty$ and taking into account the arbitrariness of $\varepsilon$, the latter relation yields by bounded convergence theorem the inequality
$$K(\lambda)\leq{\sf E}{K\left(\lambda e^{-\nu\zeta}\right)}.$$
By iteration this transforms into
\begin{equation}\label{K(lambda)leq}
K(\lambda)\leq{\sf E}{K\left(\lambda e^{-\nu Z(n)}\right)},
\end{equation}
where $Z(n):=\sum_{k=1}^{n}\zeta_k$ and $\zeta_k$, $k\in\mathbb{Z}+$, are i.i.d. random variables with the same distribution as $\zeta$.
According to strong law of large numbers and bounded convergence theorem, inequality (\ref{K(lambda)leq}) implies that
$$0\leq K(\lambda)\leq K(0+).$$
However, $K(0+)=0$ in view of (\ref{limlimK(t,lambda)=0}). Thus, Lemma~\ref{L:limK=0} is proved. $\square$

\noindent\emph{Proof of Theorem~\ref{T:main_result}} For $N=1$ and $x=w_1=0$, Theorem~\ref{T:main_result} is implied by Lemmas~\ref{L:equation}--\,\ref{L:limK=0}. Now we deal with $N>1$ and $x\in W$, say $x=w_i$. Let us discuss here the main differences between the case of single and multiple catalysts and sketch the subsequent proof omitting cumbersome details. In the multiple setting, the counterpart of integral equation (\ref{P_0(Mt_>_Lt_u)_equation}) in Lemma~\ref{L:equation} is the system of integral equations
\begin{eqnarray}
& &{\sf P}_{w_i}\left(M_t>u\right)=\alpha_i\int\nolimits_0^t{\left(1-f_i\left(1-{\sf P}_{w_i}\left(M_{t-s}>u\right)\right)\right)\,dG_i(s)}\label{NP_0(Ltgamma/Mtgamma<lambda)_equation}\\
&+&(1-\alpha_i)\sum_{j=1}^N\int\nolimits_0^t{{\sf P}_{w_j}\left(M_{t-s}>u\right)\,d\left(G_i\ast {_{W_j}\overline{F}_{w_i,w_j}(s)}\right)}+I_i^{(N)}(t;u),\nonumber
\end{eqnarray}
where $i=1,\ldots,N$ and the functions $I_i^{(N)}(t;u)$, $t\geq0$, $u\in\mathbb{R}$, have the following expression
$$(1-\alpha_i)\sum_{y\notin W}{\frac{q(w_i,y)}{q}\!\int\limits_0^t\!{\sf P}_y\left(S(t-s)>u,{_{W_k}\tau_{y,w_k}}>t-s,k=1,\ldots,N\right)\,d G_i(s)}.$$ Similar to Lemma~\ref{L:J-1(t;a)=}, each function, for $t\geq0$ and $u\geq\max\{w_1,\ldots,w_N\}$, satisfies the identity
\begin{eqnarray}\label{NJ-1(t;a)=}
\frac{q I_i^{(N)}(t;u)}{(1-\alpha_i)\beta_i}&=&{\sf P}_{w_i}\left(S(t)>u\right)-\sum_{k=1}^N\int\nolimits_0^t{{\sf P}_{w_k}\left(S(t-s)>u\right)\,d\,{_{W_k}F_{w_i,w_k}(s)}}\nonumber\\
&-&\frac{\beta_i-q}{\beta_i}\int\nolimits_0^t{{\sf P}_{w_i}\left(S(t-s)>u\right)\,d G_i(s)}\nonumber\\
&+&\sum_{k=1}^N\frac{\beta_i-q}{\beta_i}\int\nolimits_0^t{{\sf P}_{w_k}\left(S(t-s)>u\right)\,d G_i\ast{_{W_k}F_{w_i,w_k}(s)}}.
\end{eqnarray}
Now the next step is to introduce a multiple setting counterpart of function $G$ arising in Lemma~\ref{L:J-1(t;lambdaLt)ast_sim}, namely a matrix $\mathcal{G}(t)=\left(G^{N}_{i,j}(t)\right)_{i,j=1}^N$, where $G^{N}_{i,j}(t):=\delta_{i,j}\alpha_i m_i G_i(t)+(1-\alpha_i)G_i\ast{_{W_j}\overline{F}_{w_i,w_j}(t)}$, $t\geq0$, and, as usual, $\delta_{i,j}$ is the Kronecker delta. Note that the element $d_{i,j}(\lambda)$ of matrix $D(\lambda)$, $\lambda\geq0$, is just the Laplace transform of $G^{(N)}_{i,j}$.

Proceed to the multiple analogue of Lemma~\ref{L:lower_estimate} and afterwards return to the counterpart of Lemma~\ref{L:J-1(t;lambdaLt)ast_sim}. By mean value theorem on functions $f_1,\ldots,f_N$, the system of equations (\ref{NP_0(Ltgamma/Mtgamma<lambda)_equation}) implies the following vector inequality, valid coordinate-wise,
\begin{equation}\label{NP_0(Ltgamma/Mtgamma<lambda)_inequality}
\mathcal{P}(t;u)\leq\mathcal{G}\ast\mathcal{P}(t;u)+\mathcal{I}(t;u),
\end{equation}
where $\mathcal{P}(t;u):=\left({\sf P}_{w_1}\left(M_t>u\right),\ldots,{\sf P}_{w_N}\left(M_t>u\right)\right)^{\top}$ and $\mathcal{I}(t;u):=\left(I^{(N)}_1(t;u),\ldots,I^{(N)}_N(t;u)\right)^{\top}$ are the vector-columns and $\top$ stands for a matrix transposition. Recall that the operation ``$\ast$'' of convolutions of matrices is defined exactly as matrix multiplication except that we convolve elements rather than multiply them. Iterating the inequality (\ref{NP_0(Ltgamma/Mtgamma<lambda)_inequality}) $k$ times, letting $k\to\infty$ and applying Lemma~1.1 in \cite{Crump_70}, similar to (\ref{P_0(Ltgamma/Mtgamma<lambda)_inequality}) we derive
$$\mathcal{P}(t;u)\leq\sum_{k=0}^{\infty}\mathcal{G}^{\ast k}\ast\mathcal{I}(t;u).$$
Thus, as in Lemma~\ref{L:J-1(t;lambdaLt)ast_sim} for $N=1$, in case $N>1$ we inspect the asymptotic behavior of $\sum_{k=0}^{\infty}\mathcal{G}^{\ast k}\ast\mathcal{I}(t;u)$ when $u=\lambda^{-1/\gamma}L_{t+r}$ and $t\to\infty$. In full similarity to Lemma~\ref{L:J-1(t;lambdaLt)ast_sim}, employing Corollary~3.1, item (i), in \cite{Crump_70}, we deduce that
$$\left.\sum_{k=0}^{\infty}\mathcal{G}^{\ast k}\ast\mathcal{I}(t;u)\right|_{u=\lambda^{-1/\gamma}L_{t+r}}\sim \lambda e^{-\nu r}\left(K_1^{(N)},\ldots,K^{(N)}_N\right)^{\top},\quad t\to\infty.$$
The constants $K_i^{(N)}>0$, $i=1,\ldots,N$, can be written in an explicit form which is bulky and superfluous, and so omitted. Moreover, Lemma~\ref{L:lower_estimate} remains intact in case of $N>1$ as well (with, possibly, another constant $C'$ instead of $C$).

The generalization of function $J(t;u)$, $t\geq0$, $u\in\mathbb{R}$, to the case $N>1$ is the vector function $\mathcal{J}(t;u)$, $t\geq0$, $u\in\mathbb{R}$, with coordinates $J^{(N)}_i(t;u)$, $i=1,\ldots,N$, of the form
$$m_i\int\nolimits_0^t{{\sf P}_{w_i}\left(M_{t-s}>u\right)\,dG_i(s)}-\int\nolimits_0^t{\left(1-f_i\left(1-{\sf P}_{w_i}\left(M_{t-s}>u\right)\right)\right)\,dG_i(s)}.$$ A multiple setting counterpart of Lemma~\ref{L:I-1(t;lambdaLt)ast_sim} asserts that, under the same conditions, one has
$$\lim_{\lambda\to0+}\lim_{t\to\infty}\frac{1}{\lambda}\left.\sum_{k=0}^{\infty}\mathcal{G}^{\ast k}\ast\mathcal{J}(t;u)\right|_{u=\lambda^{-1/\gamma}L_t}=(0,\ldots,0)^{\top}.$$
The proof repeats that of Lemma~\ref{L:I-1(t;lambdaLt)ast_sim}, however now we apply Condition~3.1, item (i), in \cite{Crump_70} instead of Theorem~25 in \cite{Vatutin_book_09}, p.~30.

The discrepancies between both statements and proofs of Lemmas~\ref{L:limlim=theta},\ref{L:limK=0} and their counterparts are virtually insignificant. Hence we only note that the proof of the analogue of Lemma~\ref{L:limK=0} follows the one of Theorem~3.3 in \cite{Kaplan_75}. Thus, Theorem~\ref{T:main_result} is established in the case of $N\geq1$ and the starting point $x\in W$.

It remains to justify Theorem~\ref{T:main_result} in the case of $N\geq1$ and $x\notin W$. The case of starting point $x\notin W$ is  reduced to the case of $N+1$ catalysts, viz $W \,\cup \,\{x\}$, so then we can employ the results obtained for the case of $N+1$ catalysts and the starting point from $W$. Theorem~\ref{T:main_result} is proved completely. $\square$

\end{document}